# Inverting the 10-Coin Triangle Puzzle and Other Shapes: A New General Solution

Dr. Tony McCaffrey and Oscar Atwill
Eagle Hill School
Hardwick, MA

**Abstract:** This puzzle, often called the *Reverse the Triangle Puzzle*, appears regularly in puzzle books. Four rows consisting of 1 coin in row 1, 2 coins in row 2, 3 coins in row 3, and 4 coins in row 4 form the shape of a triangle. What is the fewest number of coins you need to move to flip the triangle 180 degrees so that it points in the opposite direction? In the case of 4 rows, it requires moving just 3 coins. The general solution for any number of rows can be calculated by dividing the number of coins by 3 and ignoring any remainder. A triangle of 4 rows has 10 coins, so 10/3 = 3.333 leads to 3 moves. However, this calculation is fairly nonintuitive as to why it works and requires an inelegant maneuver of periodically discarding a remainder. In this paper, we present a more intuitive formula that always calculates the answer directly and never produces a decimal answer. Further, we explore inverting other shapes (e.g., rhombus) and find an interesting connection to the formula for triangles.

## Introduction

Puzzle books often contain the *Reverse the Triangle Puzzle.* Figure 1 shows the problem with four rows, which amounts to ten coins, but puzzle books also often use the five row version. The goal is to transform the triangle on the left pointing upward into the triangle on the right pointing downward by moving the fewest number of coins possible.

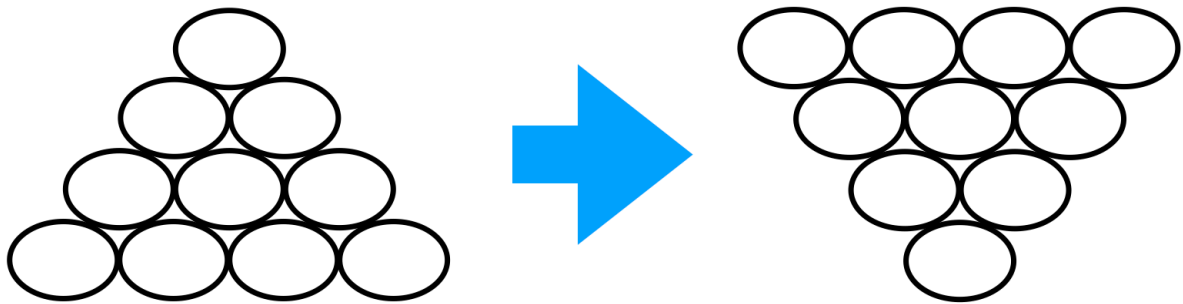

*Figure 1: Reverse the Triangle Puzzle with Four Rows*

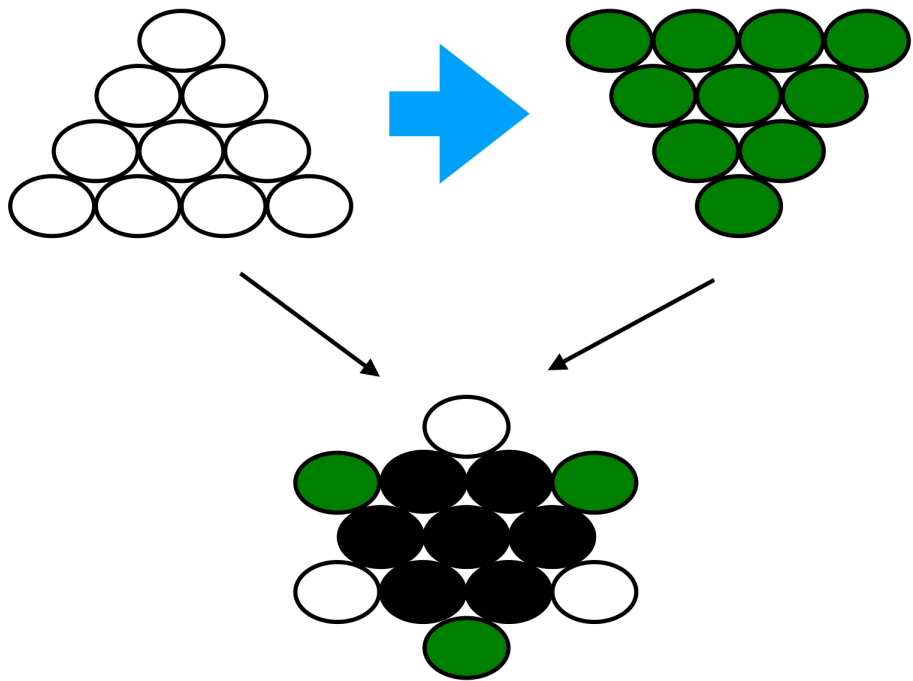

*Figure 2: Superimposing the Starting and Ending Triangles*

A great way to "see" the solution is to imagine that you superimpose the starting position of the coins and the ending position. In this way, you can easily maximize the overlap between the two and determine how many coins can remain where they are. The remaining coins must be moved.

Figure 2 shows the starting triangle as clear and the ending triangle in green. When you overlap the two triangles so that the maximum number overlap, seven coins overlap, as indicated with the color black, leaving three clear coins to be moved to flip the triangle. The three clear coins should be moved into the positions of the three green coins. This is the minimum number of coins required to flip the triangle. You can empirically try all possible overlaps to reveal that three moves is in fact the minimum.

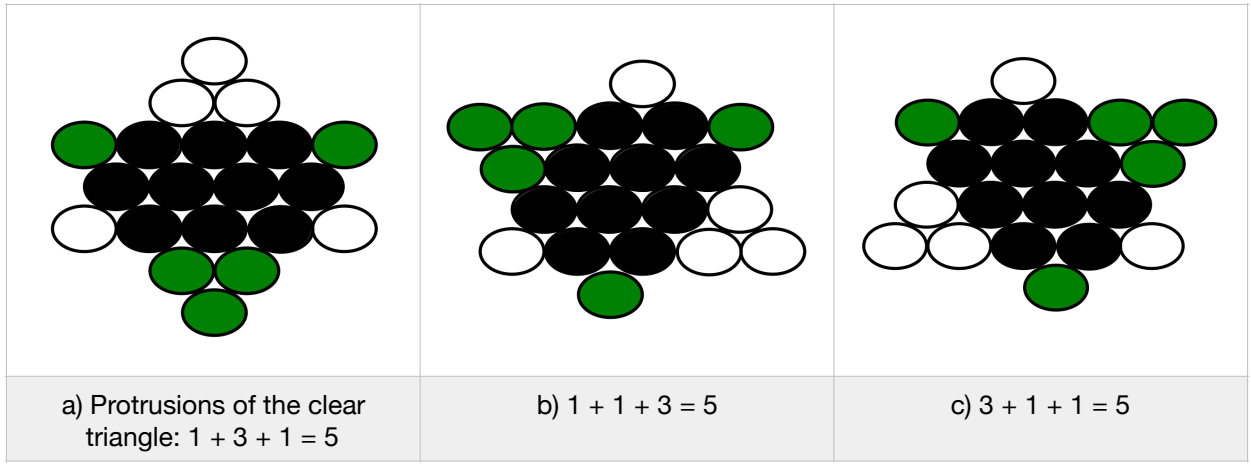

a) Protrusions of the clear triangle: 1 + 3 + 1 = 5

b) 1 + 1 + 3 = 5

c) 3 + 1 + 1 = 5

*Figure 3: Three Possible Solutions to a 5-Row Triangle*

The groups of non-overlapping coins sticking out from the center cluster of overlapping coins will be called **protrusions**. There are three clear-color protrusions sticking out from the center overlapping cluster and three green protrusions. If a protrusion consists of 1 coin, we will call it a 1-protrusion, 3 coins in a protrusion will be called a 3-protrusion, etc.

We will focus on the three clear-color protrusions since they represent the coins that need to be moved. For a five row triangle shown in Figure 3a, moving from left to right for the clear-color protrusions, there is a 1-protrusion, a 3-protrusion, and 1-protrusion. Adding them together results in 5 moves to flip the triangle. The solution in Figure 3b uses the same numbers but adds them up in a different order: 1 + 1 + 3 = 5. Figure 3c places the three numbers in yet a different order, 3 + 1 + 1 = 5. So, there may be more than one minimum solution (i.e., fewest number of moves), but each solution consists of the same size protrusions just in different orders.

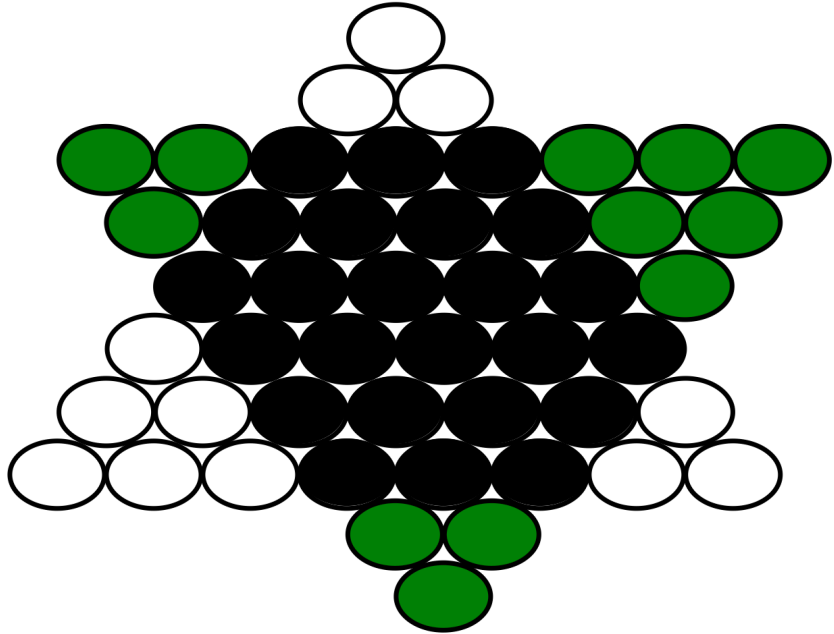

*Figure 4: One Solution to the 8-Row Triangle*

Notice that each protrusion must be a triangle of coins. Figure 3 shows 1-protrusion triangles and 3-protrusion triangles. Figure 4 shows that bigger triangles are possible—in this case a 6-protrusion triangle is required for an eight row triangle.

The number of coins in a triangular protrusion are conveniently called *triangular numbers*. They have been well-studied and the first few triangular numbers are the following: 1, 3, 6, 10, 15, 21, 28, 36, … A closed-form formula for the number of coins in a triangle of coins of $n$ rows, an $n$-protrusion, is $\frac{n(n+1)}{2}$.

Table 1: Data on Coin Triangles as Rows Increase

| Rows | Total Pennies | Old Formula | Progression of Moves | New: Sum of 3 Triangular Numbers |
|---|---|---|---|---|
| 1 | 1 | 0.3333333333 | | 0 + 0 + 0 |
| 2 | 3 | 1 | 1 - 0 = 1 | 1 + 0 + 0 |
| 3 | 6 | 2 | 2 - 1 = 1 | 1 + 1 + 0 |
| 4 | 10 | 3.3333333333 | 3 - 2 = 1 | 1 + 1 + 1 |
| 5 | 15 | 5 | 5 - 3 = 2 | 3 + 1 + 1 |
| 6 | 21 | 7 | 2 | 3 + 3 + 1 |
| 7 | 28 | 9.3333333333 | 2 | 3 + 3 + 3 |
| 8 | 36 | 12 | 3 | 6 + 3 + 3 |
| 9 | 45 | 15 | 3 | 6 + 6 + 3 |
| 10 | 55 | 18.333333333 | 3 | 6 + 6 + 6 |
| 11 | 66 | 22 | 4 | 10 + 6 + 6 |
| 12 | 78 | 26 | 4 | 10 + 10 + 6 |
| 13 | 91 | 30.333333333 | 4 | 10 + 10 + 10 |
| 14 | 105 | 35 | 5 | 15 + 10 + 10 |
| 15 | 120 | 40 | 5 | 15 + 15 + 10 |
| 16 | 136 | 45.333333333 | 5 | 15 + 15 + 15 |
| 17 | 153 | 51 | 6 | 21 + 15 + 15 |
| 18 | 171 | 57 | 6 | 21 + 21 + 15 |
| 19 | 190 | 63.333333333 | 6 | 21 + 21 + 21 |
| 20 | 210 | 70 | 7 | 28 + 21 + 21 |
| 21 | 231 | 77 | 7 | 28 + 28 + 21 |
| 22 | 253 | 84.333333333 | 7 | 28 + 28 + 28 |
| 23 | 276 | 92 | 8 | 36 + 28 + 28 |
| 24 | 300 | 100 | 8 | 36 + 36 + 28 |
| 25 | 325 | 108.33333333 | 8 | 36 + 36 + 36 |
| 26 | 351 | 117 | 9 | 45 + 36 + 36 |
| 27 | 378 | 126 | 9 | 45 + 45 + 36 |
| 28 | 406 | 135.33333333 | 9 | 45 + 45 + 45 |

Table 1, specifically Column 5, presents the three triangular numbers that must be added together to reach the number of moves required to flip triangles of various sizes. Although 0 is technically not a triangular number, for our equations we consider it to be the zeroth triangular number and a protrusion of 0 coins will be a 0-protrusion.

Column 3 of Table 1 presents the previously known solution, the calculation of division, $\frac{TotalCoins}{3}$, followed by ignoring any remainder (Baroody et al., 2004; Kindt et al., 1998). In trying to understand why this formula works, we were not completely successful but did find a few intriguing patterns. Here is one.

Column 4 presents how quickly the number of moves increase as the number of rows increase. As you can see, the number of moves increase by 1 move three times in a row, then by 2 moves three times in a row, etc. The division operation in column 3 has a remainder whenever the rate of increase is about to go up. The division has no remainder twice in a row and then on the third occurrence when the rate of increase is poised to go up, suddenly there is a remainder. We interpret this as a remainder occurring at the transition points for the rate of increase.

As you can see, we do not fully understand why Column 3 generates the correct number of moves, but we discussed one suggestive pattern that might provide a clue.

The numbers added together in Column 5 of Table 1 are all triangular numbers. We will notate the *nth* triangular number as $T_n$. The repeating pattern as we move down the fourth column in groups of three can be represented as the following. Divide the number of rows by 3. The quotient *m* indicates which triangular numbers to use (i.e., $T_m$ and/or $T_{m+1}$). The remainder *p* chooses which of the formulas to use.

$$\frac{rows}{3} = m \text{ remainder } p,$$

$$\begin{cases} 3T_m, & \text{if } p = 1 \\ T_{m+1} + 2T_m, & \text{if } p = 2 \\ 2T_{m+1} + T_m, & \text{if } p = 0 \end{cases}$$

Relying on the known fact that $T_m = \dfrac{m(m+1)}{2}$, results in the following set of three formulas.

$$\frac{rows}{3} = m \text{ remainder } p,$$

$$\begin{cases} \dfrac{3m(m+1)}{2}, & \text{if } p = 1 \\ \dfrac{(m+1)(m+2)}{2} + \dfrac{2m(m+1)}{2}, & \text{if } p = 2 \\ \dfrac{2(m+1)(m+2)}{2} + \dfrac{m(m+1)}{2}, & \text{if } p = 0 \end{cases}$$

Simplifying the algebra results in the following calculations:

$$\frac{rows}{3} = m \text{ remainder } p,$$

$$\begin{cases} \dfrac{3m^2 + 3}{2}, & \text{if } p = 1 \\ \dfrac{3m^2 + 5m + 2}{2}, & \text{if } p = 2 \\ \dfrac{3m^2 + 7m + 4}{2}, & \text{if } p = 0 \end{cases}$$

**Inverting a Rhombus**

What about shapes other than a triangle? We thought to examine a square next, but the closest we could make was a rhombus. The rhombus shape is due to how the coins are placed in the niches formed between two adjacent coins.

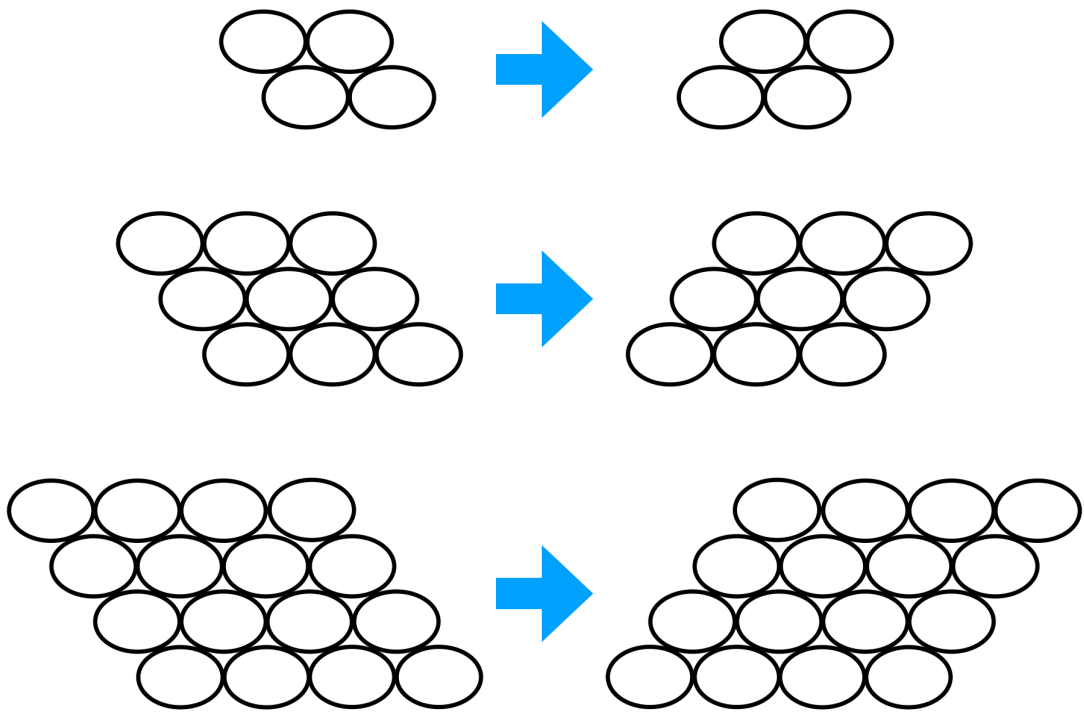

*Figure 5: Rhombi of 4, 9, and 16 coins*

With each rhombus in Figure 5, we want to change its orientation. Triangles were flipped vertically while rhombi could be flipped either vertically or horizontally and achieve the same minimum number of moves. For variety, in this paper, we will flip them horizontally. Again, we ask the same question: What is the fewest number of coins to move in order to flip a rhombus of coins horizontally?

The technique of superimposing the starting and ending positions allows us to visualize which coins need to be moved.

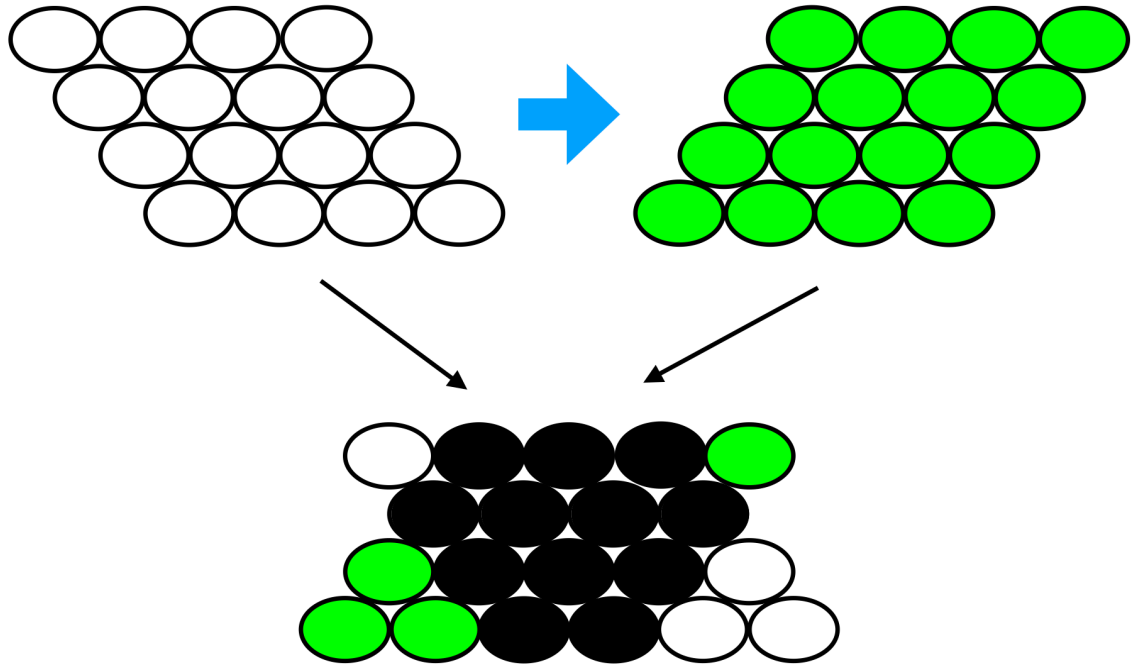

*Figure 6: Superimposing the Starting and Ending Rhombi*

Maximizing the number of overlapped coins in Figure 6 reveals that 4 coins are the fewest number required to flip the 4-row rhombus horizontally. Note that a rhombus produces only 2 protrusions that need to be moved; whereas, a triangle produced 3 protrusions. Also, each protrusion forms a triangle that consists of a triangular number of coins: 1, 3, 6, 10, 15, …

Table 2: Flipping Rhombi Horizontally

| Rows | Total Coins | Coins/4 | 2 Triangular Numbers |
|---|---|---|---|
| 1 | 1 | 0.25 | 0 + 0 |
| 2 | 4 | 1 | 1 + 0 |
| 3 | 9 | 2.25 | 1 + 1 |
| 4 | 16 | 4 | 3 + 1 |
| 5 | 25 | 6.25 | 3 + 3 |
| 6 | 36 | 9 | 6 + 3 |
| 7 | 49 | 12.25 | 6 + 6 |
| 8 | 64 | 16 | 10 + 6 |
| 9 | 81 | 20.25 | 10 + 10 |
| 10 | 100 | 25 | 15 + 10 |
| 11 | 121 | 30.25 | 15 + 15 |
| 12 | 144 | 36 | 21 + 15 |
| 13 | 169 | 42.25 | 21 + 21 |
| 14 | 196 | 49 | 28 + 21 |
| 15 | 225 | 56.25 | 28 + 28 |
| 16 | 256 | 64 | 36 + 28 |
| 17 | 289 | 72.25 | 36 + 36 |
| 18 | 324 | 81 | 45 + 36 |
| 19 | 361 | 90.25 | 45 + 45 |
| 20 | 400 | 100 | 55 + 45 |
| 21 | 441 | 110.25 | 55 + 55 |

Interestingly, if you are not bothered by discarding remainders, then $\frac{TotalCoins}{4}$, followed by ignoring any remainder will give you the answer—as shown in Table 2. With triangles, we divided by 3 and ignored the remainder to get the correct answer.

For rhombi, two triangular numbers are added together that correspond to the number of protrusions that emerge for each rhombus. As the number of rows increase in a rhombus, the pattern of adding triangular numbers repeats in cycles of two. Starting with one row (the first row of the table), the formula is $T_0 + T_0$. Two rows yields $T_1 + T_0$. Then, the cycle repeats for three and four rows, $T_1 + T_1$ and $T_2 + T_1$, respectively. The general formula is the following:

$$\frac{rows}{2} = m \text{ remainder } p,$$

$$\begin{cases} 2T_m, & \text{if } p = 1 \\ T_m + T_{m-1}, & \text{if } p = 0 \end{cases}$$

Relying on the known fact that $T_m = \frac{m(m+1)}{2}$, results in the following two formulas.

$$\frac{rows}{2} = m \text{ remainder } p,$$

$$\begin{cases} \frac{2m(m+1)}{2}, & \text{if } p = 1 \\ \frac{m(m+1)}{2} + \frac{m(m-1)}{2}, & \text{if } p = 0 \end{cases}$$

Simplifying the algebra results in the following calculations:

$$\frac{rows}{2} = m \text{ remainder } p,$$

$$\begin{cases} m^2 + m, & \text{if } p = 1 \\ m^2, & \text{if } p = 0 \end{cases}$$

**Comparing Flipping Triangles to Flipping Rhombi**

Table 3 compares the formulas for triangles to those for rhombi. Hopefully, this will help us better understand why the old formula for triangles, $\frac{TotalCoins}{3}$, works and how it compares to the comparable formula for rhombi, $\frac{TotalCoins}{4}$. In particular, since you divide by three for triangles and four for rhombi, does the divisor then necessarily represent the number of sides of the shapes involved? Curiously, for the new formulas for triangles, $\frac{rows}{3} = m$ *remainder p,* and the new formula for rhombi, $\frac{rows}{2} = m$ *remainder p,* the divisors here cannot represent the number of sides of the shape. Most likely, the divisor represents the number of protrusions that are added together.

Our work with 5-sided and 6-sided polygons of coins has not yielded patterns related to either the *Old Formula* or the *New Formula* presented in Table 3. Thus, sadly, the intriguing patterns for 3- and 4-sided polygons of coins do not generalize to higher sided polygons of coins.

Note that as the number of sides of a polygon increases, at the limit point, it approaches a circle. A circle is symmetrical from all angles, so flipping a circle will require no moves. A circle of coins in not perfectly round; but, nonetheless, it is sufficiently symmetrical from so many angles that inverting it would create the exact same shape—thus, flipping it would require no moves. So, as the number of sides increases of the polygons to be inverted, at some point they become so symmetrical from so many angles that no change in shape will take place when they are inverted. In fact, we found that the interesting formulas found for 3- and 4-sided polygons break down for 5- and 6-sided polygons of coins. We believe they will continue to break down and offer nothing interesting as the sides of the polygons of coins increase to 7-sides and higher.

Table 3: Comparison of Triangle and Rhombus Formulas

| Description | Triangle | Rhombus |
|---|---|---|
| Old Formula | $\dfrac{TotalCoins}{3}$, ignore remainder | $\dfrac{TotalCoins}{4}$, ignore remainder |
| New Formula (Part 1) | $\dfrac{rows}{3} = m$ remainder $p$ | $\dfrac{rows}{2} = m$ remainder $p$ |
| New Formula (Part 2) In terms of triangular numbers | $\begin{cases} 3T_m, & \text{if } p = 1 \\ T_{m+1} + 2T_m, & \text{if } p = 2 \\ 2T_{m+1} + T_m, & \text{if } p = 0 \end{cases}$ | $\begin{cases} 2T_m, & \text{if } p = 1 \\ T_m + T_{m-1}, & \text{if } p = 0 \end{cases}$ |
| New Formula (Part 2) Without triangular numbers | $\begin{cases} \dfrac{3m^2 + 3}{2}, & \text{if } p = 1 \\ \dfrac{3m^2 + 5m + 2}{2}, & \text{if } p = 2 \\ \dfrac{3m^2 + 7m + 4}{2}, & \text{if } p = 0 \end{cases}$ | $\begin{cases} m^2 + m, & \text{if } p = 1 \\ m^2, & \text{if } p = 0 \end{cases}$ |

**Summary**

We have created a general formula that will calculate the minimal number of moves to flip any *n*-row triangle of coins. The previous calculation for this result, divide $\dfrac{TotalCoins}{3}$ and ignore any remainder; although accurate, is not intuitive as to the pattern of behavior going on among triangular numbers as the number of rows increase. Further, the previous calculation ignores the remainder while our calculation uses the remainder to choose the appropriate formula from the three possible formulas in the pattern. In this way, we believe our new solution has merit both on intuitiveness as well as on actively using the remainder.

We also created two general formulas that will calculate the minimum moves to horizontally flip any *n*-row rhombus of coins. One formula is $\frac{TotalCoins}{4}$ and ignore any remainder. This formula regularly produces remainders that need to be ignored. The second formula adds together triangular numbers and never produces a remainder.

Finally, we examined both 5-sided and 6-sided polygons of coins and did not find any continuation of the formulas found for 3- and 4-sided coin polygons. Further, we reasoned that no interesting formulas will be present for 7-sided polygons of coins and higher.